# A New Fast Direct Method For Solving Quasi-Toeplitz System of Equations

PREPRINT

**Shahin Hasanbeigi**
Department of Applied Mathematics
Tarbiat Modares University
Tehran, Iran
shahin77hb@gmail.com

**Abstract**

The objective of this study is to present a novel, efficient, and fast direct method for solving linear systems of equations whose coefficient matrix is a tridiagonal Quasi-Toeplitz matrix. Such matrices are frequently encountered in the discretization of second-order differential equation problems with Neumann boundary conditions, the discretization of Quasi-Birth-and-Death processes known as QBD matrix equations, and other related applications. The presented method has been demonstrated to produce favorable results in terms of run time for moderately large systems when compared to classical direct methods, such as the LU and PLU methods.

## 1 Introduction

Consider the following nonsingular linear system of equation

$$Tx = b \tag{1}$$

where T is an $n \times n$ Quasi-Toeplitz matrix of the form

$$T = \begin{bmatrix} \lambda & \beta & & & & & \\ \gamma & \alpha & \beta & & & & \\ & \gamma & \alpha & \beta & & & \\ & & \ddots & \ddots & \ddots & & \\ & & & \ddots & \ddots & \ddots & \\ & & & & \gamma & \alpha & \beta \\ & & & & & \gamma & \delta \end{bmatrix} \tag{2}$$

Toeplitz matrices and their variants, including Quasi-Toeplitz and Banded Toeplitz, represent a distinct category of structured matrices that find applications in diverse fields such as applied mathematics, engineering, and scientific computing. The unique structure of Toeplitz matrices has been the subject of extensive research and publications in the domain of linear system solvers, see [1, 2].

In general, there are two approaches for solving Toeplitz linear systems: direct methods and iterative methods. Iterative methods are primarily comprised of classical splitting iteration methods and Krylov subspace iteration methods, such as GMRES, MINRES, and so forth. Conversely, direct methods, including LU and PLU methods, as discussed in reference [2], are frequently applicable to small and moderate-sized problems. However, they are often too costly to be practical for large sparse problems. Additionally, direct methods may exhibit numerical instability and a loss of accuracy in the solutions they produce [3].

Tridiagonal Toeplitz matrices are a common matrix form that arises in a variety of scientific and engineering applications. As an illustration, the discretization of second-order differential equations will result in the formation of a linear system, wherein the coefficient matrix takes on the form of a tridiagonal Toeplitz or Quasi-Toeplitz matrix [4]. The special structure of Toeplitz and Quasi-Toeplitz matrices allows for the development of highly efficient algorithms that are significantly faster than classical algorithms such as LU and PLU methods. Some of these types of algorithms can be found in [5, 6].

This paper is organized as follows. The following section introduces the new method and its implementation for solving the tridigonal Quasi-Toeplitz system (1). The final section compares the results of the new method for solving several tridigonal quasi-Toeplitz systems, including those derived from the discretization of second-order differential equations with Neumann boundary conditions, with the LU and PLU methods.

## 2 New Direct Method

In this section we first assume that $n$ is a even number. The new method does not change under this assumption, but when $n$ is odd there will be a few small diffeences in some steps that will be discused in the final part of this section.

So let $n$ be even and,

$$J = \begin{bmatrix} & & & & & 1 \\ 1 & & & & & \\ & 1 & & & & \\ & & \ddots & & & \\ & & & \ddots & & \\ & & & & 1 & \\ & & & & & 1 \end{bmatrix}$$

It's very easy to see that $J$ is a orthonormal matrix, that mean

$$JJ^T = J^T J = I$$

So we have

$$JTJ^T Jx = Jb \tag{3}$$

We can rename matrix and vectors in (3) as follow

$$T^{(1)} x^{(1)} = b^{(1)} \tag{4}$$

Where

$$T^{(1)} = JTJ^T = \left[\begin{array}{cc|cccccc} \delta_1 & 0 & 0 & & \cdots & 0 & \gamma \\ 0 & \lambda_1 & \beta & 0 & \cdots & & 0 \\ \hline 0 & \gamma & \alpha & \beta & & & \\ & 0 & \gamma & \alpha & \beta & & \\ \vdots & \vdots & & \ddots & \ddots & \ddots & \\ 0 & & & & \gamma & \alpha & \beta \\ \beta & 0 & & & & \gamma & \alpha \end{array}\right] \equiv \left[\begin{array}{c|c} D^{(1)} & P^{(1)} \\ \hline W^{(1)} & A^{(1)} \end{array}\right] \quad , \quad \lambda_1 = \lambda \ , \ \ \delta_1 = \delta \tag{5}$$

Considering $2 \times 2$ block LU factorization of matrix $T^{(1)}$, we have

$$T^{(1)} = \left[\begin{array}{c|c} D^{(1)} & P^{(1)} \\ \hline W^{(1)} & A^{(1)} \end{array}\right] = \begin{bmatrix} I_2 & \\ W^{(1)} D^{(1)^{-1}} & I_{n-2} \end{bmatrix} \begin{bmatrix} D^{(1)} & P^{(1)} \\ & A^{(1)} - W^{(1)} D^{(1)^{-1}} P^{(1)} \end{bmatrix} \tag{6}$$

Also we consider below partitioning for $x^{(1)}$ and $b^{(1)}$:

$$x^{(1)} = Jx = \begin{bmatrix} x_1^{(1)} \\ \\ x_2^{(1)} \end{bmatrix} \quad , \quad b^{(1)} = Jb = \begin{bmatrix} b_1^{(1)} \\ \\ b_2^{(1)} \end{bmatrix}$$

where $x_1^{(1)}$ , $x_2^{(1)}$ , $b_1^{(1)}$ and $b_2^{(1)}$ have propper dimentions corespond to $T^{(1)}$ patitioning in (5).

By replacing above partitioning and block factorization of $T^{(1)}$ in (3), we reach to the below equations:

$$\begin{cases} D^{(1)}x_1^{(1)} + P^{(1)}x_2^{(1)} & = b_1^{(1)} \\ \big(A^{(1)} - W^{(1)}D^{(1)^{-1}}P^{(1)}\big)x_2^{(1)} & = b_2^{(1)} - W^{(1)}D^{(1)^{-1}}b_1^{(1)} \end{cases} \tag{7}$$

Second equation in (7) is a $(n-2)\times(n-2)$ system of equations that should be solved for $x_2^{(1)}$. After computing $x_2^{(1)}$ from second equation, $x_1^{(1)}$ can be computed easily with first equation by multiplying $D^{(1)^{-1}}$ to it from the left. Notice that computing $D^{(1)^{-1}}$ is very cheap and easy since $D^{(1)}$ is a small $2\times 2$ diagonal matrix and its inverse simply is

$$D^{(1)^{-1}} = \begin{bmatrix} \frac{1}{\delta_1} & 0 \\ 0 & \frac{1}{\lambda_1} \end{bmatrix}$$

As stated above, computing $x_2^{(1)}$ require solving a $(n-2)\times(n-2)$ system of equations that its coefficient matrix is

$$A^{(1)} - W^{(1)}D^{(1)^{-1}}P^{(1)} = \begin{bmatrix} \lambda_2 & \beta & & & & & \\ \gamma & \alpha & \beta & & & & \\ & \gamma & \alpha & \beta & & & \\ & & \ddots & \ddots & \ddots & & \\ & & & \ddots & \ddots & \ddots & \\ & & & & \gamma & \alpha & \beta \\ & & & & & \gamma & \delta_2 \end{bmatrix} \tag{8}$$

where $\lambda_2 = \alpha - \frac{\beta\gamma}{\lambda_1}$ and $\delta_2 = \alpha - \frac{\beta\gamma}{\delta_1}$ . It is clear that above matrix is a Quasi-Toeplitz matrix like (2), therefore solving second equation in (7) is similar to solving the original system but its dimention is $n-2$ and it's smaller than the original system. So, we multiply matrix $J$ with approprate size to the second equation similar to (3):

$$J\big(A^{(1)} - W^{(1)}D^{(1)^{-1}}P^{(1)}\big)J^T Jx_2^{(1)} = J\big(b_2^{(1)} - W^{(1)}D^{(1)^{-1}}b_1^{(1)}\big) \tag{9}$$

Again, with renaming matrix ans vectors we have:

$$T^{(2)}x^{(2)} = b^{(2)} \tag{10}$$

where

$$T^{(2)} = J\big(A^{(1)} - W^{(1)}D^{(1)^{-1}}P^{(1)}\big)J^T = \left[\begin{array}{cc|cccccc} \delta_2 & 0 & 0 & & \cdots & & 0 & \gamma \\ 0 & \lambda_2 & \beta & 0 & \cdots & & & 0 \\ \hline 0 & \gamma & \alpha & \beta & & & & \\ & 0 & \gamma & \alpha & \beta & & & \\ \vdots & \vdots & & \ddots & \ddots & \ddots & & \\ 0 & & & & \gamma & \alpha & \beta & \\ \beta & 0 & & & & \gamma & \alpha & \end{array}\right] \equiv \left[\begin{array}{c|c} D^{(2)} & P^{(2)} \\ \hline W^{(2)} & A^{(2)} \end{array}\right]$$

and just like before, we consider below partitioning for $x^{(2)}$ and $b^{(2)}$:

$$x^{(2)} = Jx_2^{(1)} = \begin{bmatrix} x_1^{(2)} \\ \\ x_2^{(2)} \end{bmatrix} \quad , \quad b^{(2)} = J\big(b_2^{(1)} - W^{(1)}D^{(1)^{-1}}b_1^{(1)}\big) = \begin{bmatrix} b_1^{(2)} \\ \\ b_2^{(2)} \end{bmatrix} \tag{11}$$

Similar to $T^{(1)}$, we consider $2\times 2$ block LU factorization of matrix $T^{(2)}$

$$T^{(2)} = \left[\begin{array}{c|c} D^{(2)} & P^{(2)} \\ \hline W^{(2)} & A^{(2)} \end{array}\right] = \begin{bmatrix} I_2 & \\ W^{(2)}D^{(2)^{-1}} & I_{n-4} \end{bmatrix} \begin{bmatrix} D^{(2)} & P^{(2)} \\ & A^{(2)} - W^{(2)}D^{(2)^{-1}}P^{(2)} \end{bmatrix} \tag{12}$$

and by replacing it into (11) we will have

$$\begin{cases} D^{(2)}x_1^{(2)} + P^{(2)}x_2^{(2)} & = b_1^{(2)} \\ \big(A^{(2)} - W^{(2)}D^{(2)^{-1}}P^{(2)}\big)x_2^{(2)} & = b_2^{(2)} - W^{(2)}D^{(2)^{-1}}b_1^{(2)} \end{cases} \tag{13}$$

Exactly like explained before, if we compute $x_2^{(2)}$ from second equation of (13), we can easily compute $x_1^{(2)}$ from the first equation since

$$D^{(2)^{-1}} = \begin{bmatrix} \frac{1}{\delta_2} & 0 \\ 0 & \frac{1}{\lambda_2} \end{bmatrix}$$

Computing $x_2^{(2)}$ require solving a $(n-4) \times (n-4)$ system of equations that its coefficient matrix is

$$A^{(2)} - W^{(2)} D^{(2)^{-1}} P^{(2)} = \begin{bmatrix} \lambda_3 & \beta & & & & & \\ \gamma & \alpha & \beta & & & & \\ & \gamma & \alpha & \beta & & & \\ & & \ddots & \ddots & \ddots & & \\ & & & \ddots & \ddots & \ddots & \\ & & & & \gamma & \alpha & \beta \\ & & & & & \gamma & \delta_3 \end{bmatrix} \quad , \ \lambda_3 = \alpha - \frac{\beta\gamma}{\lambda_2 \lambda_1} \ , \ \delta_3 = \alpha - \frac{\beta\gamma}{\delta_2 \delta_1} \tag{14}$$

which is Quasi-Toeplitz martix like the original system that should be solved.

So far, we can see from above precedure that every time we repeat the above steps, size of the system that has computational expenses to solve goes down by 2. Also every time we repeat above steps, we need to solve another $2 \times 2$ system involving $D^{(i)}$'s which are not computationally expensive and very quick since $D^{(i)}$'s are $2 \times 2$ diagonal mtrices and we do not need further computations for $D^{(i)^{-1}}$'s.

If we repeat above steps $t = \frac{n-4}{2} + 1$ times we will have:

$$T^{(t)} x^{(t)} = b^{(t)} \tag{15}$$

where

$$T^{(t)} = \left[\begin{array}{cc|cc} \delta_t & 0 & 0 & \gamma \\ 0 & \lambda_t & \beta & 0 \\ \hline 0 & \gamma & \alpha & \beta \\ \beta & 0 & \gamma & \alpha \end{array}\right] = \left[\begin{array}{c|c} D^{(t)} & P^{(t)} \\ \hline W^{(t)} & A^{(t)} \end{array}\right] \tag{16}$$

and

$$\lambda_t = \alpha - \frac{\beta\gamma}{\lambda_{t-1}} \ , \ \delta_t = \alpha - \frac{\beta\gamma}{\delta_{t-1}}$$

Again, we consider $2 \times 2$ block LU factorization of $T^{(t)}$ and by replacing it into (15) we have:

$$\begin{cases} D^{(t)} x_1^{(t)} + P^{(t)} x_2^{(t)} & = b_1^{(t)} \\ \left(A^{(t)} - W^{(t)} D^{(t)^{-1}} P^{(t)}\right) x_2^{(t)} & = b_2^{(t)} - W^{(t)} D^{(t)^{-1}} b_1^{(t)} \end{cases} \tag{17}$$

Second equation of (17) is $2 \times 2$ system which solving it gives us $x_2^{(t)}$. By replacing $x_2^{(t)}$ into first equation in (17) we can compute $x_1^{(t)}$ and by puting them together we have

$$x^{(t)} = \begin{bmatrix} x_1^{(t)} \\ \\ x_2^{(t)} \end{bmatrix}$$

Notice that, according to (11), we named $x^{(i)}$'s as follow:

$$x^{(i+1)} = J x_2^{(i)}$$

where matrix $J$ has approprate size according to $x_2^{(i)}$. So according to this, after computing $x^{(t)}$ from (17), we have

$$x_2^{(t-1)} = J^T x^{(t)}$$

Now, since we have $x_2^{(t-1)}$, we use equation

$$D^{(t-1)}x_1^{(t-1)} + P^{(t-1)}x_2^{(t-1)} = b_1^{(t-1)}$$

and compute $x_1^{(t-1)}$, and by putting them together we acuire

$$x^{(t-1)} = \begin{bmatrix} x_1^{(t-1)} \\ \\ x_2^{(t-1)} \end{bmatrix}$$

We can repeat above steps until we compute $x^{(1)}$. Then, from (4), we have

$$x = J^T x^{(1)}$$

which is the solution of our original system of equations (1).

In the mentioned method, at each step, we proceed to the second subsystem of the previous system and continue this process until we reach the $4 \times 4$ system in (16). Then we solve the system (16) according to the (17). Once the solution to (16) has been obtained, with the help of matrix $J$ and the first subsystem in every step, we compute the solution of previous step and continue this until we reach to the solution of our original system. One important thing that we should notice in this precedure is that in every step, besides the size of matrices and vectors, also some elements of them will change. For example, in case of $W^{(i)}$'s and $P^{(i)}$'s, position of $\beta$ and $\gamma$ will remain the same bot the size of these matrices will decrease until they reach to their size in (16). But in case of the right hand side vectors $b^{(i)}$, in every step the first and last elements of it will change. When our original coefficient matrix is a general tridiagonal matrix, following these changes is exceedingly challenging, However, given that the original coefficient matrix is a Quasi-Toeplitz matrix, it is possible to leverage its structure and derive explicit formulas for the dimensions and elements of matrices and vectors at each stage of the method.

At the beginning of this section we put $\lambda_1 = \lambda$ , $\delta_1 = \delta$ and we saw that

$$\lambda_2 = \alpha - \frac{\beta\gamma}{\lambda_1} \quad , \quad \delta_2 = \alpha - \frac{\beta\gamma}{\delta_1}$$

It's easy to verify that for $i = 2, ..., t$:

$$\lambda_i = \alpha - \frac{\beta\gamma}{\lambda_{i-1}} \quad , \quad \delta_i = \alpha - \frac{\beta\gamma}{\delta_{i-1}} \tag{18}$$

Furthermore, for $i = 1, 2, ..., t$:

$$D^{(i)} = \begin{bmatrix} \delta_i & 0 \\ 0 & \lambda_i \end{bmatrix}_{2\times 2} \tag{19}$$

Also, for $i = 1, 2, ..., t$:

$$P^{(i)} = \begin{bmatrix} 0 & 0 & \cdots & 0 & \gamma \\ \beta & 0 & \cdots & 0 & 0 \end{bmatrix}_{2\times n-2i} \quad , \quad W^{(i)} = \begin{bmatrix} 0 & \gamma \\ 0 & 0 \\ \vdots & \vdots \\ 0 & 0 \\ \beta & 0 \end{bmatrix}_{n-2i\times 2} \quad , \quad A^{(i)} = \begin{bmatrix} \alpha & \beta & & & \\ \gamma & \alpha & \beta & & \\ & \ddots & \ddots & \ddots & \\ & & \gamma & \alpha & \beta \\ & & & \gamma & \alpha \end{bmatrix}_{n-2i\times n-2i} \tag{20}$$

By looking at (9) and (10) it is clear that for $i = 2, ..., t$:

$$b^{(i)} = J\big(b_2^{(i-1)} - W^{(i-1)} {D^{(i-1)}}^{-1} b_1^{(i-1)}\big)$$

so it is obvious that $\beta, \gamma$ , $\delta_i$'s and $\lambda_i$'s, are effective on elements of $b^{(i)}$'s. For $i = 1$ we have:

$$b_1^{(1)} = \begin{bmatrix} b_n \\ b_1 \end{bmatrix} \quad , \quad b_2^{(1)} = \begin{bmatrix} b_2 \\ \vdots \\ b_{n-1} \end{bmatrix} \tag{21}$$

For $i = 2, ..., t$, it is not very hard to confirm that:

$$b_1^{(i)} = \begin{bmatrix} b_{n-(i-1)} & + & \sum\limits_{k=1}^{i-1}(-1)^k \dfrac{\beta^k b_{n-(i-1)+k}}{\prod\limits_{j=1}^{k}\delta_{i-j}} \\ b_i & + & \sum\limits_{k=1}^{i-1}(-1)^k \dfrac{\gamma^k b_{i-k}}{\prod\limits_{j=1}^{k}\lambda_{i-j}} \end{bmatrix} \quad , \quad b_2^{(i)} = \begin{bmatrix} b_{i+1} \\ \vdots \\ b_{n-i} \end{bmatrix} \tag{22}$$

With help of above equations, we can compress our method in an algorithm as follow:

**Algorithm 1** Tridiagonal Quasi-Toeplitz Solver *(when n is a even number)*

**Require:** $\alpha, \beta, \gamma, \lambda, \delta, \mathbf{b}$

Let: $t = \frac{n-4}{2} + 1$ , $\lambda_1 = \lambda$ , $\delta_1 = \delta$ , $b_1^{(1)} = \begin{bmatrix} b_n \\ b_1 \end{bmatrix}$ , $b_2^{(1)} = \begin{bmatrix} b_2 \\ \vdots \\ b_{n-1} \end{bmatrix}$

**for** $i = 2, \cdots, t$ **do**

  $\lambda_i = \frac{\beta\gamma}{\lambda_{i-1}}$ , $\delta_i = \frac{\beta\gamma}{\delta_{i-1}}$ , $b_1^{(i)} = \begin{bmatrix} b_{n-(i-1)} & + & \sum\limits_{k=1}^{i-1}(-1)^k \dfrac{\beta^k b_{n-(i-1)+k}}{\prod\limits_{j=1}^{k}\delta_{i-j}} \\ b_i & + & \sum\limits_{k=1}^{i-1}(-1)^k \dfrac{\gamma^k b_{i-k}}{\prod\limits_{j=1}^{k}\lambda_{i-j}} \end{bmatrix}$ , $b_2^{(i)} = \begin{bmatrix} b_{i+1} \\ \vdots \\ b_{n-i} \end{bmatrix}$

**end for**

Let : $A^{(t)} = \begin{bmatrix} \alpha & \beta \\ \gamma & \alpha \end{bmatrix}$

Determine $W^{(t)}$ and $P^{(t)}$ by (20), and $D^{(t)}$ by (19)

Solve: $\left(A^{(t)} - W^{(t)}D^{(t)^{-1}}P^{(t)}\right)x_2^{(t)} = b_2^{(t)} - W^{(t)}D^{(t)^{-1}}b_1^{(t)}$ for $x_2^{(t)}$.

**for** $i = 1, \cdots, t$ **do**

  $x_1^{(t-i+1)} = D^{(t-i+1)^{-1}}\left(b_1^{(t-i+1)} - P^{(t-i+1)}x_2^{(t-i+1)}\right)$

  $x^{(t-i+1)} = \begin{bmatrix} x_1^{(t-i+1)} \\ x_2^{(t-i+1)} \end{bmatrix}$

  $x_2^{(t-i)} = J^T x^{(t-i+1)}$

**end for**

**return** $x = x_2^{(0)}$

Algorithm (1) is designed to solve tridiagonal Quasi-Toeplitz systems that its dimentions are an even number. When dimention of a system is an odd number, the generallity of new method does not change, but a few small modifications needed to adjust the algorithm.

In case of $n$ being odd, we need to repeat our precedure $t = \frac{n-3}{2} + 1$ times, and we will have:

$$T^{(t)}x^{(t)} = b^{(t)} \quad , \quad T^{(t)} = \left[\begin{array}{cc|c} \delta_t & 0 & \gamma \\ 0 & \lambda_t & \beta \\ \hline \beta & \gamma & \alpha \end{array}\right] = \left[\begin{array}{c|c} D^{(t)} & P^{(t)} \\ \hline W^{(t)} & A^{(t)} \end{array}\right] \tag{23}$$

Using $2 \times 2$ block LU factorization of above matrix we will have:

$$\begin{cases} D^{(t)}x_1^{(t)} + P^{(t)}x_2^{(t)} & = b_1^{(t)} \\ \left(A^{(t)} - W^{(t)}D^{(t)^{-1}}P^{(t)}\right)x_2^{(t)} & = b_2^{(t)} - W^{(t)}D^{(t)^{-1}}b_1^{(t)} \end{cases} \tag{24}$$

The difference between (24) and (17) is that, the second equation of (24) is a scaler equation, thus

$$x_2^{(t)} = \frac{b_2^{(t)} - W^{(t)}D^{(t)^{-1}}b_1^{(t)}}{A^{(t)} - W^{(t)}D^{(t)^{-1}}P^{(t)}} \tag{25}$$

After computing $x_2^{(t)}$ by (25), backsubstitution precedure that gives us solution $x$ is exactly according to the algorithm(1). Other differences are in $W^{(t)}$, $P^{(t)}$ and $A^{(t)}$ matrices. In every step, dimentions of these matrices will reduce by 2, same as algorithm(1), but in the last step these matrices will be as following:

$$P^{(t)} = \begin{bmatrix} \gamma \\ \beta \end{bmatrix} \ , \ W^{(t)} = \begin{bmatrix} \beta & \gamma \end{bmatrix} \ , \ A^{(t)} = \alpha$$

Therefore, algorithm for our new method when dimention of our system is odd, is following algorithm,

---

**Algorithm 2** Tridiagonal Quasi-Toeplitz Solver *(when n is a odd number)*

---

**Require:** $\alpha, \beta, \gamma, \lambda, \delta, \mathbf{b}$

Let: $t = \frac{n-3}{2} + 1 \ , \ \lambda_1 = \lambda \ , \ \delta_1 = \delta \ , \ b_1^{(1)} = \begin{bmatrix} b_n \\ b_1 \end{bmatrix} \ , \ b_2^{(1)} = \begin{bmatrix} b_2 \\ \vdots \\ b_{n-1} \end{bmatrix}$

**for** $i = 2, \cdots, t$ **do**

  $\lambda_i = \frac{\beta\gamma}{\lambda_{i-1}} \ , \ \delta_i = \frac{\beta\gamma}{\delta_{i-1}} \ , \ b_1^{(i)} = \begin{bmatrix} b_{n-(i-1)} + \sum_{k=1}^{i-1} (-1)^k \frac{\beta^k b_{n-(i-1)+k}}{\prod_{j=1}^{k} \delta_{i-j}} \\ b_i + \sum_{k=1}^{i-1} (-1)^k \frac{\gamma^k b_{i-k}}{\prod_{j=1}^{k} \lambda_{i-j}} \end{bmatrix} \ , \ b_2^{(i)} = \begin{bmatrix} b_{i+1} \\ \vdots \\ b_{n-i} \end{bmatrix}$

**end for**

Let : $A^{(t)} = \alpha \ , W^{(t)} = \begin{bmatrix} \beta & \gamma \end{bmatrix} \ , \ P^{(t)} = \begin{bmatrix} \gamma \\ \beta \end{bmatrix}$ and determine $D^{(t)}$'s by (19).

Compute $x_2^{(t)}$ with (25).

**for** $i = 1, \cdots, t$ **do**

  $x_1^{(t-i+1)} = {D^{(t-i+1)}}^{-1} \left( b_1^{(t-i+1)} - P^{(t-i+1)} x_2^{(t-i+1)} \right)$

  $x^{(t-i+1)} = \begin{bmatrix} x_1^{(t-i+1)} \\ x_2^{(t-i+1)} \end{bmatrix}$

  $x_2^{(t-i)} = J^T x^{(t-i+1)}$

**end for**

**return** $x = x_2^{(0)}$

---

So in general, algorithm for solving tridiagonal Quasi-Toeplitz system can be express as follow:

---

**Algorithm 3** Tridiagonal Quasi-Toeplitz Solver

---

**Require:** $\alpha, \beta, \gamma, \lambda, \delta, \mathbf{b}, \mathbf{n}$

**if n** is even **then**

  Use Algorithm(1).

**else**

  Use Algorithm(2).

**end if**

---

One important notice that can be inferred from our new algorithm, is that new method works perfectly for tridiagonal Toeplitz systems. The term Quasi-Toeplitz is not a necessary condition for the new method and it is simply an extention of Toeplitz matrices. Applying new method to a tridiagonal Toeplitz system only differs in the first step, and from second step, coeffiient matrix will turn into a Quasi-Toeplitz matrix.

On the other hand, most of the times, our tridiagonal Quasi-Toeplitz systems are results of discretizing different types of problems, for example second order differential equation with Neumann boundary condition. Therefore, different kinds of discretizing methods can end with different types of matrix equation. Even if we use same discretizing method, such as Finite Difference method, the final tridiagonal Quasi-Toeplitz system may be different. For example, by using finite difference method for discretizing a second order differential equation with Neumann boundary condition, some might reach to the system

$$Tx = b$$

where

$$T = \begin{bmatrix} \alpha & \lambda & & & & & \\ \gamma & \alpha & \beta & & & & \\ & \gamma & \alpha & \beta & & & \\ & & \ddots & \ddots & \ddots & & \\ & & & \ddots & \ddots & \ddots & \\ & & & & \gamma & \alpha & \beta \\ & & & & & \delta & \alpha \end{bmatrix}$$

Since our new method is designed for coefficient matrices like (2), by multiplying matrix

$$Z = \begin{bmatrix} \beta/\lambda & & & & \\ & 1 & & & \\ & & \ddots & & \\ & & & 1 & \\ & & & & \gamma/\delta \end{bmatrix}$$

same as [6], we will have:

$$\hat{T} = ZT = \begin{bmatrix} \frac{\alpha\beta}{\lambda} & \beta & & & & & \\ \gamma & \alpha & \beta & & & & \\ & \gamma & \alpha & \beta & & & \\ & & \ddots & \ddots & \ddots & & \\ & & & \ddots & \ddots & \ddots & \\ & & & & \gamma & \alpha & \beta \\ & & & & & \gamma & \frac{\alpha\gamma}{\delta} \end{bmatrix}$$

which is a tridiagonal Quasi-Toeplitz matrix like (2). Therefore, even if our system is not a tridiagonal Quasi-Toeplitz like (2), we can use a matrix similar to $Z$ to reform the original system before using the new method.

It is also worth noting that the new method is more efficient when the coefficient matrix is diagonally dominant. As with all direct methods, the new method will be affected by error emission more rapidly when $T$ is not a diagonally dominant matrix, particularly as the size of $T$ increases. When T is diagonally dominant, the new method will maintain the precision of the solution to a high degree of accuracy, comparable to that of other classical methods such as LU factorization. However, the new method will also compute the solution at a significantly faster rate. The following section will present a series of examples to illustrate this point.

## 3 Numerical Experiments

In this section we use some examples to show effectiveness of our new algorithm.

All the numerical tests were done on an ASUS laptop PC with AMD A12 CPU, 8Gb RAM and by Matlab R2016(b) with a machine precision of $10^{-16}$. For convenience, throughout our numerical experiments, we denote the relative residual error of methods by Relative error $= ||b - Ax||/||b||$, and computing time of methods by Time (in seconds). In all examples, LU denotes the LU facorization method with pivoting (if its necessary). In all tables, the Time is the average value of computing times required by performing the corresponding algorithm 10 times.

**Example 1 :**

In our first example we use some artificial systems to demonstrate effectiveness of new algorithm when the coefficient matrix $T$ is diagonally dominant and when is not diagonally dominant. In this example we persume $b = rand(n, 1)$. Also, dimention of systems is been choosen such that $n$ be even and odd, and in a increasing order.

First, let

$$\alpha = 1\ ,\ \ \beta = 2\ ,\ \ \gamma = 3\ ,\ \ \lambda = 4\ ,\ \ \delta = 5$$

Table 1: Results for different $n$ when $T$ is not diagonally dominant

| n | New method Time | New method Relative error | LU Time | LU Relative error |
|---|---|---|---|---|
| $2^5$ | 0.0034 | 1.5732 e -14 | 0.0031 | 2.9117 e -14 |
| 97 | 0.0060 | 1.7225 e -9 | 0.0179 | 1.5908 e -9 |
| $2^7$ | 0.0107 | 3.1350 e -7 | 0.0307 | 1.5745 e -6 |
| 183 | 0.0210 | 2.7394 e -2 | 0.0845 | 4.7652 e -2 |
| $2^{10}$ | 5.2837 | very big number | 30.6960 | very big number |

Table (1) shows the results for above system. Since $T$ is not diagonally dominant, as the size of the matrix increases, methods start loosing their percision. Table (1) results shows that new method's percision lost is very similar to PLU's percision lost. But, almost every time, new method's time is better than PLU, even when size of $T$ is large and methods are not working.

For next test let $T$ be a diagonallly dominant matrix. Let

$$\alpha = -4\ ,\quad \beta = 1\ ,\quad \gamma = 1\ ,\quad \lambda = 2\ ,\quad \delta = 3$$

Table 2: Results for different $n$ when $T$ is diagonally dominant

| n | New method Time | New method Relative error | LU Time | LU Relative error |
|---|---|---|---|---|
| $2^5$ | 0.0037 | 2.1439 e -16 | 0.0031 | 1.4987 e -16 |
| 97 | 0.0049 | 1.9532 e -16 | 0.0196 | 1.7044 e -16 |
| $2^7$ | 0.0119 | 1.9602 e -16 | 0.0294 | 1.5908 e -16 |
| 183 | 0.0155 | 1.7225 e -16 | 0.0832 | 1.6009 e -16 |
| $2^{10}$ | 7.0742 | 1.9000 e -16 | 45.4517 | 1.6401 e -16 |

Table(2) results clearly shows that when $T$ is diagonally domonant, methods keep their percision, even when size of the matrix is fairly large. Results also shows that, almost every time, new method's running time is way lower than LU method. So as been mentioned before, new method shows best run time in compair to other direct methods, while its percision stays same as them. Therefore, new method is more suitable for solving tridiagonal Quasi-Toeplitz systems, or tridiagonal Toeplitz systems, in compair to other direct methods.

**Example 2 :**

In this example we consider a second order differential equation with Neumann boundary condition.

$$\begin{cases} u''(x) - u(x) = 0 \\ u'(0) = -1\ ,\ u'(1) = \frac{-1}{e} \end{cases} \quad ,\ \Omega = [0,1]$$

Exact solution of above equation is $u(x) = e^{-x}$. With use of derivative approximations

$$u''(x_i) = \frac{u(x_{i+1}) - 2u(x_i) + u(x_{i-1})}{h^2}\ ,\ u'(x_i) = \frac{u(x_{i+1}) - u(x_i)}{h}\ ,\ u'(x_i) = \frac{u(x_i) - u(x_{i-1})}{h}$$

we would have:

$$\begin{cases} u_2 - u_1 = -h & i = 1 \\ u_{i+1} + \alpha u_i + u_{i-1} = 0 & i = 2, \cdots, N-1 \\ u_N - u_{N-1} = \frac{-h}{e} & i = N \end{cases}$$

where $h$ is stepsize, $N = \frac{1}{h}$, and $\alpha = -(2 + h^2)$. Above recursive formula give us the system

$$\begin{bmatrix} -1 & 1 & & & & & \\ 1 & \alpha & 1 & & & & \\ & 1 & \alpha & 1 & & & \\ & & \ddots & \ddots & \ddots & & \\ & & & \ddots & \ddots & \ddots & \\ & & & & 1 & \alpha & 1 \\ & & & & & 1 & -1 \end{bmatrix} \begin{bmatrix} u_1 \\ u_2 \\ \\ \vdots \\ \\ u_{N-1} \\ u_N \end{bmatrix} = \begin{bmatrix} -h \\ 0 \\ \\ \vdots \\ \\ 0 \\ \frac{h}{e} \end{bmatrix}$$

Table(3) shows the results for solving example 2 by new method and LU factorization.

Table 3: Results for example 2

| h | New method Time | New method Relative error | LU Time | LU Relative error |
|---|---|---|---|---|
| 0.2 | 0.0024 | 3.2326 e -15 | 8.8434 e -4 | 1.3282 e -15 |
| 0.1 | 9.9007 e -4 | 2.5255 e -15 | 0.0014 | 3.4361 e -15 |
| 0.05 | 9.1278 e -4 | 6.4756 e -15 | 0.0018 | 5.2745 e -15 |
| 0.02 | 0.0024 | 2.9712 e -14 | 0.0057 | 2.8807 e -14 |
| 0.01 | 0.0052 | 9.0089 e -14 | 0.0202 | 9.6402 e -14 |

**Example 3 :**

For our last example, consider below differential equation:

$$\frac{d^2y}{dx^2} - \frac{T}{EI}y = \frac{qx(L-x)}{2EI} \quad , \quad y'(0) = 0 \quad , \quad y'(L) = 0$$

This differential equation governs the deflection $y$ in a simply supported beam with a uniform load $q$ and a tensile axial load $T$ where

$x$ = location along the beam (in)

$T$ = tension applied (lbs)

$E$ = Youngs modulus of elasticity of the beam (psi)

$I$ = second moment of area (in4 )

$q$ = uniform loading intensity (lb/in)

$L$ = length of beam (in)

given

$$T = 7200\ lbs\ ,\ q = 5400\ lbs/in\ ,\ L = 75\ in\ ,\ E = 30\ Msi\ ,\ \text{and}\ I = 120\ in^4$$

the goal is to solve the differential equation for $y$. With use of derivative approximations

$$y''(x_i) = \frac{y(x_{i+1}) - 2y(x_i) + y(x_{i-1})}{h^2} \ ,\ y'(x_i) = \frac{y(x_{i+1}) - y(x_i)}{h} \ ,\ y'(x_i) = \frac{y(x_i) - y(x_{i-1})}{h}$$

we would have:

$$\begin{cases} y_2 - y_1 = 0 & i = 1 \\ y_{i+1} + \alpha y_i + y_{i-1} = \frac{h^2 x_i (L-x_i)}{2EI} & i = 2, \cdots, N-1 \\ y_N - y_{N-1} = 0 & i = N \end{cases}$$

where $h$ is stepsize, $N = \frac{L}{h}$, and $\alpha = -(2 + \frac{h^2 T}{EI})$. Above recursive formula give us the system

$$\begin{bmatrix} -1 & 1 & & & & & \\ 1 & \alpha & 1 & & & & \\ & 1 & \alpha & 1 & & & \\ & & \ddots & \ddots & \ddots & & \\ & & & \ddots & \ddots & \ddots & \\ & & & & 1 & \alpha & 1 \\ & & & & & 1 & -1 \end{bmatrix} \begin{bmatrix} y_1 \\ y_2 \\ \\ \vdots \\ \\ y_{N-1} \\ y_N \end{bmatrix} = \begin{bmatrix} 0 \\ \frac{h^2 x_2 (L-x_2)}{2EI} \\ \\ \vdots \\ \\ \frac{h^2 x_{N-1} (L-x_{N-1})}{2EI} \\ 0 \end{bmatrix}$$

Table(4) shows the results for solving example 3 by new method and LU factorization.

Table 4: Results for example 3

| h | New method Time | New method Relative error | LU Time | LU Relative error |
|---|---|---|---|---|
| 25 | 0.0028 | 4.9685 e -16 | 0.0013 | 0 |
| 15 | 0.0030 | 2.0086 e -16 | 0.0014 | 8.9827 e -17 |
| 5 | 0.0033 | 1.4034 e -16 | 0.0019 | 4.0512 e -17 |
| 1 | 0.0073 | 1.5583 e -16 | 0.0109 | 8.4302 e -17 |
| 0.5 | 0.0146 | 1.4552 e -16 | 0.0519 | 1.0000 e -16 |
| 0.1 | 2.2500 | 5.7751 e -15 | 7.2511 | 3.4779 e -15 |

# 4 Conclusion

In this paper, we introduced a new fast direct method for solving tridigonal Quasi-Toeplitz systems.these types of matrix equations are common in discretizing second order differential equation with Neumann boundary condition, QBD matrix equation, which is stem from discretized Quasi-Birth-and-Death processes. Additionally, they have applications in various scientific and engineering fields. Furthermore, it is notable that QBD matrix equations frequently manifest as block tridiagonal Quasi-Toeplitz systems, which represent a prospective focus of my future research endeavors.